\documentclass[11pt, twoside]{article}
\usepackage{amssymb, amsmath, enumerate, latexsym, theorem}
\usepackage{graphicx,url}

\newcommand{\Flux}{\operatorname{Flux}}
\newcommand{\Hel}{\operatorname{Hel}}
\newcommand{\curl}{\operatorname{curl}}
\newcommand{\Twist}{\operatorname{Twist}}
\newcommand{\Tw}{\operatorname{Tw}}
\newcommand{\pf}{\operatorname{pf}}
\newcommand{\bt}{\mathbf{t}}
\newcommand{\cR}{\mathcal R}
\newcommand{\GR}{\mathcal{G}_R}
\newcommand{\D}{\mathcal{D}}
\newcommand{\Z}{\mathbf{Z}}
\newcommand{\R}{\mathbf{R}}

\newtheorem{theo}{Theorem}
\newtheorem{prop}{Proposition}[section]

\newtheorem{lemma}[prop]{Lemma}
\newtheorem{question}[prop]{Question}

\theoremstyle{definition}
\newtheorem{example}[prop]{Example}
\AtEndEnvironment{example}{\null\hfill$\diamond$}%
\newtheorem{defin}[prop]{Definition}

\AtEndEnvironment{remark}{\null\hfill$\diamond$}%

\title{Domino tilings beyond 2D}

\author{Caroline J. Klivans and Nicolau Saldanha}

\begin{document}

\maketitle

\begin{abstract}
There is a rich history of domino tilings in two dimensions.
Through a variety of techniques we
can answer questions such as:
how many tilings are there of a given region or what does the space of all
tilings look like?
These questions and their answers become significantly more difficult in
dimension three and above.
Despite this curse of
dimensionality, there have been exciting recent advances in the theory.
Here we briefly review foundational results of two-dimensional domino tilings
and highlight where challenges arise when generalizing to higher dimensions.
We then survey results in  higher dimensional domino tilings,
focusing on the question of connectivity of the space of tilings, 
and the open problems that still remain.
\end{abstract}

\section{Introduction}

There is a rich literature on two dimensional domino tilings.  The theory has been studied from a variety of perspectives including combinatorics, probability and statistical physics.  
In two dimensions, we consider covering 
quadriculated 
regions in the plane by $2 \times 1$ tiles consisting of two
squares
glued together. Equivalently, we seek a  perfect matching in planar subgraphs of the $\Z^2$ lattice.
Fundamental questions in the area include: Which regions can be tiled? Given a tileable region, how many tilings are there? and Is there a natural way to move between one tiling and another?  These and other questions have been investigated and answered in many ways using various interesting techniques.

In three or more dimensions,
we consider covering
cubiculated
regions of $\R^d$ by domino tiles consisting of two
unit
cubes glued together.
Equivalently, we seek a perfect matching in subgraphs of $\Z^d$. 
The seemingly subtle  move to higher dimensions or equivalently higher degree lattices proves to add great complexity.
While most of the questions from two dimensions translate naturally, many of the techniques used to answer them are no longer sufficient.
In higher dimensions, we have very few answers to any of the questions above.  Nonetheless, much interesting progress has been made recently, especially in the development of new tools and techniques.  Here we survey the recent advancements in higher dimensional tilings, providing many partial answers and stating the open problems that remain.

We start by briefly reviewing two dimensions and then
introduce three dimensional domino tiling. In three dimensions, it is
already more challenging to simply have visual representations of tilings.  We provide examples and include data for small regions in order to appreciate
the size and complexity of tilings that appears already in these small cases.  The article highlights the question of connectedness of the space of tilings; i.e. given a region can we move from one tiling to the other via a sequence of local moves?  A topological statistic called the twist proves to be quite relevant for this question and others.  We see the twist throughout the survey and again include computations in order to demonstrate the difficulty of the higher dimensional theory.  We will also touch on the topics of random tilings, enumeration, and moving beyond three dimensions.  

\bigbreak

The reader will notice that we give a lot of space to our own work.  This bias is mostly caused by our own interests and expertise.  There is some very interesting work about which we did not write in great detail, but we hope that the material here might lead an interested reader to investigate more into the area of higher dimensional tiling.

{The first author acknowledges support from the NSF.
 The second author acknowledges support from CNPq, CAPES and Faperj.}

\section{Domino tilings in two dimensions}

A \textit{quadriculated region} $\cR$ is a finite union of unit squares
with disjoint interiors and fitting together as the squares
$[i,i+1] \times [j,j+1]$ in the plane.
A \textit{domino} is the union along an edge
of two adjacent unit squares of $\cR$.
A \textit{domino tiling} is a set of dominoes
with disjoint interiors
which covers the region $\cR$.

Equivalently, we may construct a dual graph $\cR^\ast$.
Vertices of $\cR^\ast$ are centers of unit squares of $\cR$
and edges of $\cR^\ast$ join centers of unit squares of $\cR$
which are adjacent along an edge of $\cR$.
Thus, $\cR^\ast$ is a bipartite graph with several nice properties.
A domino tiling then corresponds to a perfect matching of $\cR^\ast$;
this is also called a \textit{dimer} cover.

\begin{figure}[ht]
\centerline{\includegraphics[width=3in]{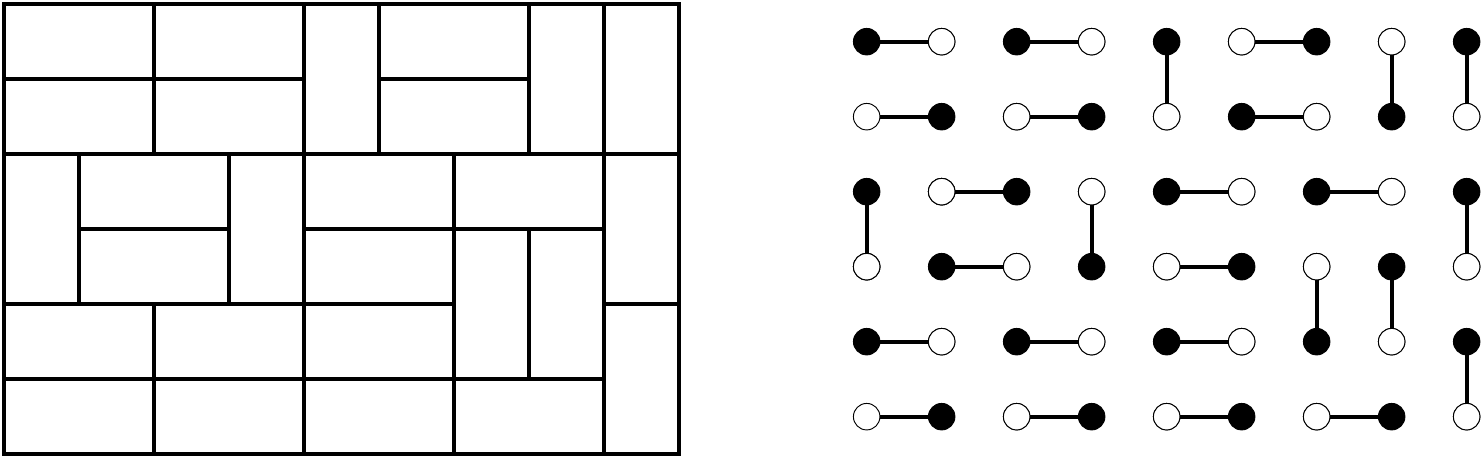}}
\caption{A two dimensional tiling and
the corresponding dimer cover of the dual graph.}  
\label{fig:2DTiling}
\end{figure}

The theory of domino tilings in two dimensions spans many disciplines
and draws on various techniques including from combinatorics,
probability, and statistical physics.
We provide here only a very quick overview of this vast subject and refer the reader to the article by Krattenthaler in this same volume and surveys~\cite{kenyon_2004, MR1370508}.

A natural first question is enumerative. 

\begin{question}
Let $\cR$ be a planar quadriculated region.
How many domino tilings of $\cR$ are there?
\end{question}

This question was addressed in Kasteleyn~\cite{KASTELEYN19611209}
and also in Temperley and Fisher~\cite{MR136398}.
Their methods rely on linear algebra:
Pfaffians, determinants and the transfer matrix-method.
More explicitly, given a region one constructs its adjacency matrix
$A$ (with entries in $\{0,1\})$.
There is a correspondence between tilings of $\cR$ and
monomials in the expansion of $\det A$
but there is also unwanted cancelling.
By changing signs of entries in a clever way
 a matrix $K$ is obtained (with entries in $\{0,\pm 1\}$)
for which $\det K$ gives the number of tilings of $\cR$.

\begin{figure}[ht]
   \centerline{\includegraphics[width=3in]{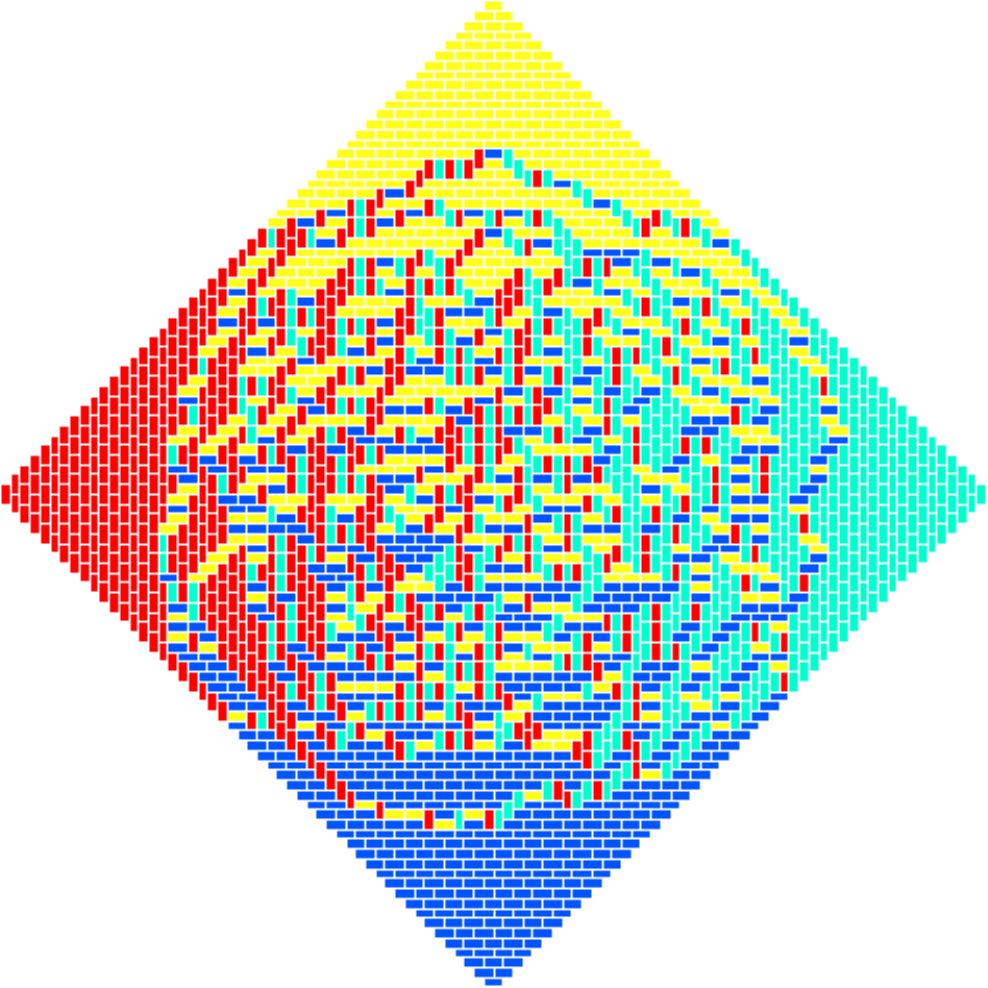}}
   \caption{The Aztec diamond: a random tiling evinces the Arctic Circle}
   \label{fig:arctic}
\end{figure}

For some special regions the methods yield useful formulas.
For instance,
the number $N(m,n)$ of domino tilings of the $m \times n$ rectangle is given by the somewhat surprising expression:
$$ N(m,n) = 
\prod_{j=1}^{\lceil{\frac{m}{2}}\rceil}
\prod_{k=1}^{{\lceil{\frac{n}{2}}\rceil}}
\left(4\cos^2 \frac{\pi j}{(m+1)} + 4\cos^2\frac{\pi k}{n+1} \right).$$
This formula can be used to obtain both estimates of $N(m,n)$
and arithmetical properties of this number.
Other regions would be studied later,
such as the Aztec diamond, shown in Figure~\ref{fig:arctic}.

\begin{question} What does a random tiling look like?
\end{question}

This is a deep question and it has received a lot of attention.
We only mention a few notable contributions.

The \textit{Aztec diamond} is the region shown in Figure~\ref{fig:arctic}.
The Arctic circle theorem by {Jockush, Propp and Shor~\cite{Aztec1}} and the work~\cite{Aztec2} by Cohn, Elkies and Propp 
shows that there are (at least) two different kinds of behavior for tiles in such regions.
Inside of a circle inscribed in the diamond
the tiling is disordered (as a liquid or a gas);
outside of the same circle, the tiling is frozen
into a solid wall of dominoes, all with the same orientation.

This beautiful result inspired a lot of research.
For instance: a similar behavior occurs in other regions,
as discussed in~\cite{Amoebae} by Kenyon, Okounkov  and Sheffield.

Another natural question is the following.  

\begin{question}
Given a region $\cR$, can we move from one tiling to another?
Is the space of tilings connected via local moves? 
\end{question}

\begin{figure}[ht]
    \centerline{
    \includegraphics[width=0.25\linewidth]{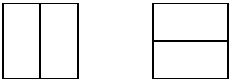}}
    \caption{A flip move.}
    \label{fig:flipmove}
\end{figure}

We first define the simplest kind of local move: a \textit{flip}.
In order to perform a flip,
remove two adjacent parallel dominoes
(forming a $2 \times 2$ block)
and place them back in the only other possible manner
(rotated by $\pi/2$), see Figure~\ref{fig:flipmove}.

A fundamental result here is the following, shown by Thurtson~\cite{MR1072815}.

\begin{theo}
In two dimensions,
any two tilings of a simply connected region are flip connected.
\end{theo}

Similar results hold for regions which are not simply connected,
such as tori or planar regions with holes, but other properties must then be taken into account.

The above questions and others have inspired much beautiful mathematics.  In this article, we will pose the same questions for domino tilings in higher dimensions.  Generally speaking we have far fewer answers in higher dimensions, but will discuss recent progress in various directions.

\section{Basics of three-dimensional domino tilings}
\label{sec:basics}

Our survey will primarily focus on domino tilings in three dimensions.
The regions of consideration are cubical polyhedral complexes $\cR$,
usually embedded in $\R^3$, safely thought of as collections of unit cubes
which intersect along an entire face, edge or corner, or not at all.
The dual graph $\mathcal{G}(\cR^*)$ of $\cR$ is formed by taking a vertex
for each cube in $\cR$ and an edge for any pair of cubes sharing a side.
The graph $\mathcal{G}(\cR^*)$ is a subgraph of the $\Z^3$ lattice and induces a bipartite two coloring of the cubes in $\cR$ -- color a cube white if the sum of the coordinates of the corresponding
vertex is even and black otherwise.

We will assume certain niceness conditions on our regions requiring
that $\cR$ is a topological manifold and that $\cR$ is balanced, i.e. the
number of white cubes is equal to the number of black cubes, otherwise a domino tiling would not be possible. 

A three-dimensional domino is a $2 \times 1 \times 1$ union of two
adjacent cubes in $\cR$.  A three-dimensional tiling of $\cR$ is a
collection of domino tilings which exactly cover $\cR$.  Two dominos may
only intersect along faces of their boundary and the union of all
dominos must equal $\cR$.  A domino in $\cR$ corresponds to an edge in
$\mathcal{G}(\cR^*)$.  A domino tiling of $\cR$ corresponds to a perfect
matching of $\mathcal{G}(\cR^*)$, sometimes called a dimer cover of $\mathcal{G}(\cR^*)$.

\begin{figure}[ht]
\centerline{
\includegraphics[width=1.5in]{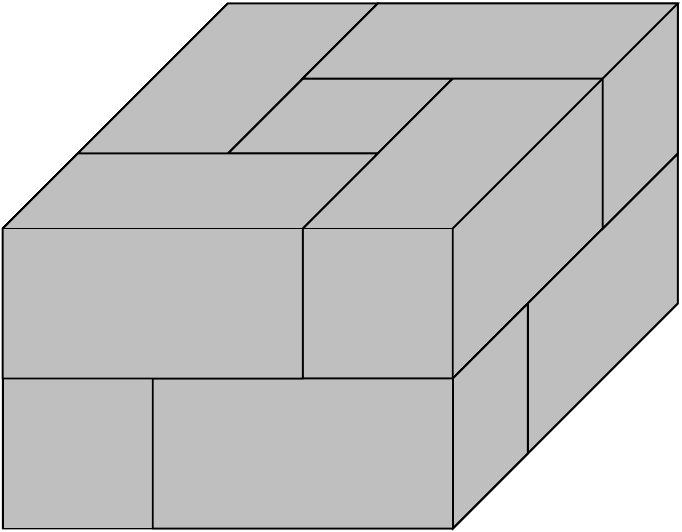} 
\hspace{.7in}
\includegraphics[width=1.5in]{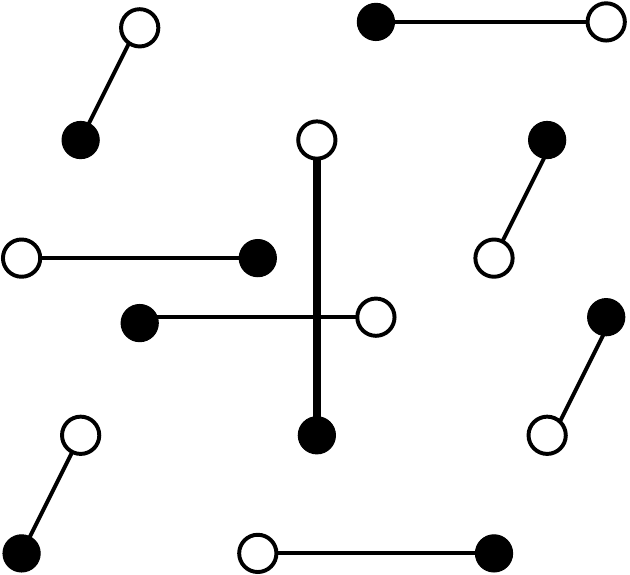}}
 \caption{Three dimensional dominoes are  $2 \times 1 \times 1$  bricks of two adjacent cubes.  On the left is a three dimensional tiling of the $3 \times 3 \times 2$ box.  On the right is the corresponding dimer cover of the dual graph.}
\label{fig:3Dregion}
\end{figure}

\begin{figure}[ht]
      \centerline{    \includegraphics[width=12cm]{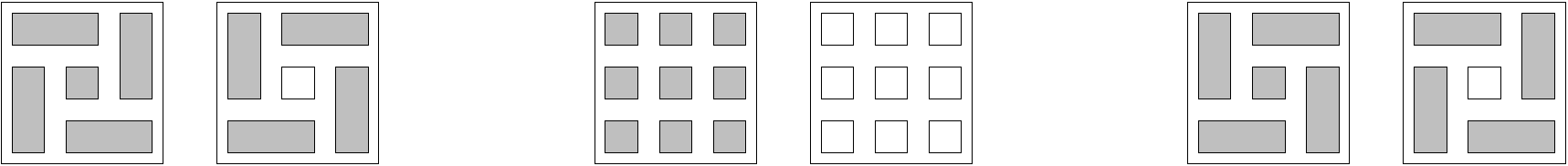}}
   \caption{Representations of tilings of the $3 \times 3\times 2$ box  by floors.  }
 \label{fig:floors}
    \end{figure}

As we can see already in this small example, it is hard to visualize
domino tilings in three dimensions.  To better visualize such tilings,
we use a system of floors, as illustrated in Figure~\ref{fig:floors}.
The region $\cR$ is broken into slices or floors by fixing values of the
$z$-direction.
The left picture of Figure~\ref{fig:floors} shows the
{same tiling we saw in Figure~\ref{fig:3Dregion},
with floors parallel to the $xy$ plane.}
Domino tiles that sit
properly in a single plane are drawn within that plane.  Vertical
domino tiles that span two planes are drawn as two corresponding
squares, one black, one white. 
The image in the middle
{of Figure~\ref{fig:floors}}
is the domino tiling of $\cR$ with
all vertical dominos.
{The third tiling is obtained by reflection.}

\subsection{Questions}

As we saw above, for a simply connected region in two dimensions, any
two tilings are connected by a series of flip moves.  The flip move is
natural in three dimensions as well.  A \emph{flip} removes two
adjacent parallel dominoes and places them back rotated within a $2
\times 2 \times 1$ block.  

We are led to the question: Given a region in three
dimensions, are any two tilings connected by a series of flip moves?
Or equivalently, in three dimensions is the space of tilings connected
by flips?  Unlike two dimensions, the answer is no.  Of course the two
dimensional case required the region to be simply connected. In three
dimensions, even for nice regions, the answer is still no.  For
example,
{the first and last tilings in Figure~\ref{fig:floors}}
both have no possible flip moves.

\begin{figure}[ht]
\centerline{\includegraphics[width=12cm]{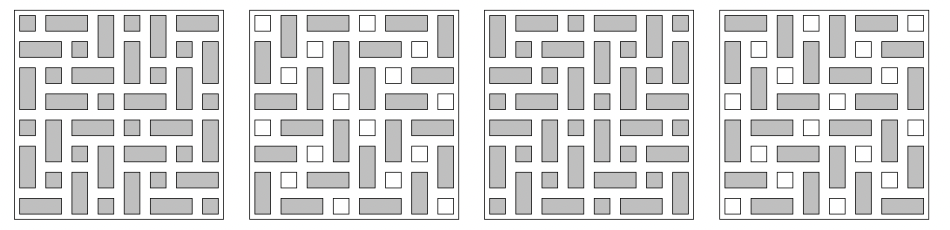}}
\caption{A tiling of the $8 \times 8 \times 4$ box
with no possible flip moves.}
\label{fig:Noflip4}
\end{figure}

The $3 \times 3 \times 2$ box is a small region
but it is not the size of the box that is causing the obstruction.
There are also examples of tilings of larger boxes that similarly have
no flip moves.
For even $N$, if $L$ and $M$ are multiples of $4$,
it is easy to adapt the example in Figure~\ref{fig:Noflip4}
to obtain tilings of the $L \times M \times N$ box
which admit no flips.
{For large boxes,
it would be interesting to know more
about the set of tilings 
which admit no flips.}

In fact, the flip connected components of nice regions
can be quite complicated. 
In the table below, we show the number and size
of the flip connected components for the 
{ $3 \times 3 \times 2$ box} and the       
{ $4 \times 4 \times 4$ box}. \\

\bigbreak

\centerline{\large{The number of flip connected components of small boxes}}
\begin{tabular}{l l | l l}
  \hline \hline
  & { $3 \times 3 \times 2$ box} & &  { $4 \times 4 \times 4$ box} \\
  \hline
& Number of tilings: $229$ & &  Number of tilings: $5,051,532,105$\\
  & Connected components: $3$ & &  Connected components: $93$\\
  \hline
  & Sizes: & & Sizes: \\
  & 227 & & 4,412,646,453 \\
  & 1 & & 2 $\times$ 310,185,960\\
  & 1 & & 2 $\times$ 8,237,514\\
  & & &      2 $\times$ 718,308\\
  & & &       2 $\times$ 283,044 \\
  & & &       6 $\times$ 2,576\\
  & & &       24 $\times$ 618\\
  & & &       $24 \times 236$\\
  & & &        6 $\times$ 4\\
  & & &       $24 \times 1$ \\
 \end{tabular} \\

For the $3 \times 3 \times 2$ box, the two components of size one are
shown in Figure~\ref{fig:floors},
{first and third tiling.}
There are an additional $227$
tilings of this size box, all of which are connected by a series of
flip moves. 
For the $4 \times 4 \times 4$ box, components are listed by size and multiplicity.  For example, there are $6$ components of size $2,576$ and $24$ components of size $1$; i.e. tilings with no flip moves. 
         All of this leads to the following motivating question:
         \begin{question}
           \label{qn:local}
  For three dimensional regions,  how and when can we move from one tiling to another via local moves?
\end{question}
Although Question~\ref{qn:local} is not yet fully resolved,
there has been significant progress recently
which we outline in Section~\ref{sec:local}.

The data in the table also reveals that there are a very large number of tilings of
even quite small regions.  This leads to another natural question.  

\begin{question}
  \label{qn:count}
For a fixed region $\cR$ in three (or more) dimensions, how many tilings are there of $\cR$? 
\end{question}

 Question~\ref{qn:count} is also not resolved; we discuss the enumerative perspective of higher dimensional domino tiling in Section~\ref{sec:count}.

\section{Connectivity by local moves}
\label{sec:local}

As we have seen above, unlike in the two-dimensional case, the
collection of three dimensional tilings of a fixed region is not
generally connected by flip moves.  The flip move is, in essence,
two-dimensional.  For three dimensional regions, a new move
is needed.

         \subsection{A Three Dimensional Move}

\begin{figure}[ht]
   \centerline{\includegraphics[width=1.65in]{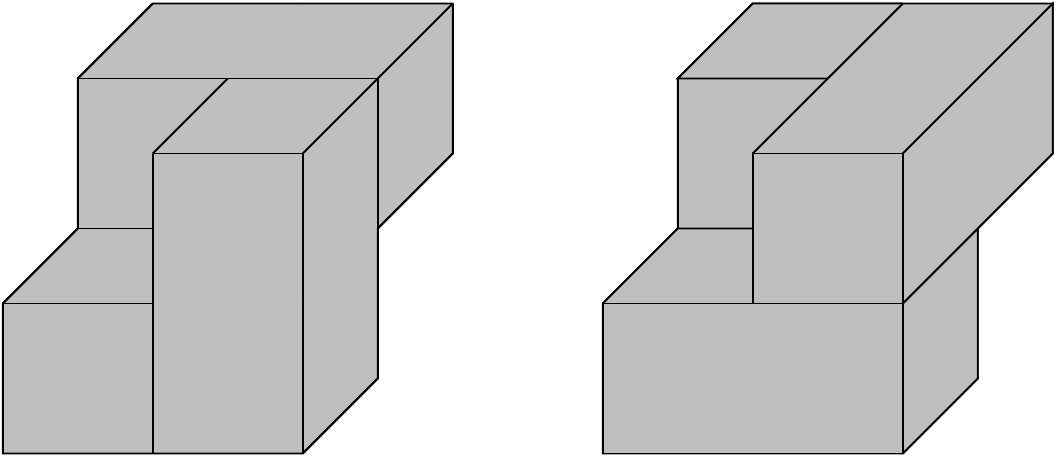}}
\caption{A positive trit move replaces the $3$ dominos seen on the left with the three dominos on the right. A negative trit moves replaces the dominos on the right with those on the left.}
 \label{fig:trit}
  \end{figure}

  \begin{defin}
A \emph{trit} move replaces $3$ dominoes, one parallel to each axis inside a $2 \times 2 \times 2$ box, see Figure~\ref{fig:trit}.  
  \end{defin}
  
  The trit move is fundamentally three dimensional in nature. It
  alters a domino in each of the three axis directions.  Moreover, it cannot be reproduced by a sequence of flip moves.  Instead, it can often be used as a local move for a tiling with no flip moves.
For example, the
{first and third} tilings
of Figure~\ref{fig:floors}
have no flip moves but can be connected via a sequence of trit and flip moves.  The driving open question in this area is whether or not these two local moves, flips and trits, suffice to connect three dimensional regions.

\begin{question}
\label{qn:tritflip}
{Let $\cR$ be a fixed
three-dimensional region.
Assume that $\cR$ is connected and simply connected.
Is the set of all tilings of $\cR$
connected by flips and trits?}
\end{question}

\begin{figure}[ht]
\centerline{ \includegraphics[width=12cm]{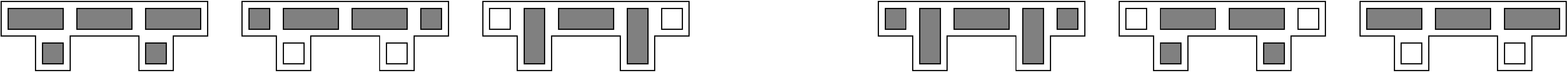}}
\caption{Two tilings of a very special region.}
\label{fig:counter}
\end{figure}

In general, the answer is no.
{Figure \ref{fig:counter} shows a counterexample:
no flip or trit is possible.}
Notice that if the region is not simply connected then 
the flux (to be defined below) is an invariant under both flips and trits.

However, the problem may be that these regions are not ``nice'' enough. 
In two dimensions,
regions need to be simply connected
to guarantee that you can flip between all possible tilings.  
In three dimensions,
it is not clear what the correct notion of a nice region is. 
The cylinders mentioned in the previous paragraph are contractible,
so the notion can not be of a purely topological nature.
Whatever the notion, we assume it includes boxes. 
Even for these nicest of regions, the question of connectedness remains open.
 
\begin{question}
\label{qn:localbox}
   Are all tilings of the  $[0,L] \times [0,M] \times [0,N]$ box
connected by flips and trits?
\end{question}


\subsection{A Topological Approach}

Our first approach to
{studying such questions}
is topological, using homology, knot theory
and a cubical approximation of surfaces. 

For this approach we define two topological invariants. 
The \emph{Flux} of a tiling can be thought of as  ``flow across surfaces''
and takes on values in $H_1(\cR;\mathbf{Z})$. 
The \emph{Twist} of a tiling can be thought of as how much a tiling
is ``twisted by trits''. 
Twist is an integer if Flux = $0$,
otherwise the Twist is in $\mathbf{Z} / m \mathbf{Z}$  where
{the value of}
$m$ depends on Flux.  
We will define these more precisely shortly.  

For this approach, we also need the notion of \emph{refinement}. 
A \emph{refinement of a region} decomposes each cube
into $5 \times 5 \times 5$ smaller cubes. 
A \emph{refinement of a tiling} decomposes each domino
into $5 \times 5 \times 5$ smaller dominoes,  each parallel to the original. 
Refinement gives us a little extra space to perform local moves. 
First we note the following. 

\begin{prop}
If $t_0$ and $t_1$ are two tilings of a region connected by 
flips (rsp.  flips and trits) then
their refinements are also connected by 
flips (rsp.  flips and trits).
\end{prop}

The converse is false, examples exist in
{the $4 \times 4 \times 2$ or}
$4 \times 4 \times 4$ boxes.

Let $\cR$ be a three dimensional region and
$t_0, t_1$ be two tilings of $\cR$. 
The \emph{difference of tilings} is defined as follows:

 $$t_1 - t_0 := \textrm{union of tiles (with the orientation of $t_0$ reversed)}.$$

\begin{figure}[ht]
  \centerline{ \includegraphics[width=1.0in]{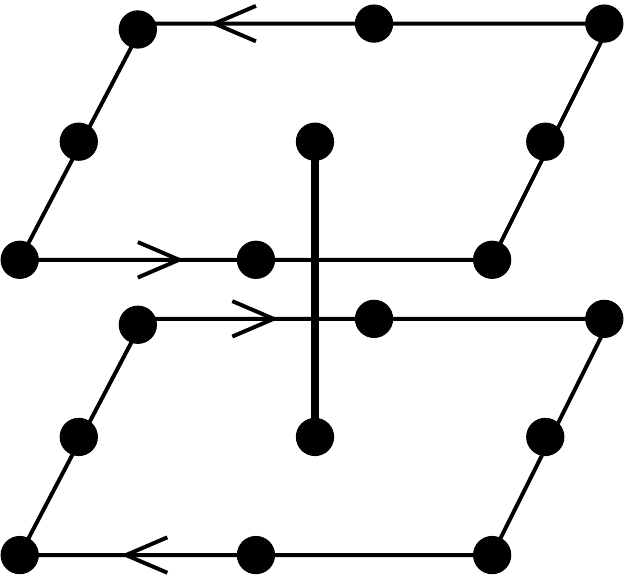}}
\caption{The system of cycles formed from the difference of tilings seen in Figure~\ref{fig:floors}.}
\label{fig:cycles}
   \end{figure}

The difference of two tilings yields a system of cycles in the dual graph $\mathcal{G}(\cR^*)$, we ignore trivial $2$-cycles formed when a tile is the same in both tilings, see Figure~\ref{fig:cycles}.  Homologically, the difference forms a boundary, $t_1 - t_0 \in Z_1(\cR^*; \mathbf{Z})$.

We are now ready to define the flux,
the first of our two topological invariants. 
Given a region $\cR$, 
fix a base tiling $t_{\oplus}$. 
For the box, we often choose the tiling
with all vertical tiles for $t_{\oplus}$. 
The flux of a tiling $t$ is the homology class
of the difference between $t$ and the base tiling. 

\begin{defin} For a region $\cR$ with base tiling $t_{\oplus}$,
the flux of a tiling $t$ is 
$$\Flux(t) := [t - t_{\oplus}] \in H_1(\cR^*; \mathbf{Z}).$$
\end{defin}

In order to see flux as a concept of flow and in order to define
the second topological invariant, we borrow a notion from knot theory. 
In knot theory, a Seifert surface is a
{connected and oriented compact} surface
whose boundary is a given knot or link.
It is a classical theorem that any link has a Seifert surface associated to it.
In our setting, the difficulty is that we are working with
discrete spaces and not smooth curves. 
Our boundaries are differences of tilings and
the surfaces we want must live in $\cR^*$.  

\begin{defin}
A (discrete) Seifert surface for a pair $t_0, t_1$
is a connected embedded oriented topological surface $S$
with boundary  $t_0 - t_1$. 
\end{defin}

\begin{centering}
\begin{figure}[ht]
\centerline{ \includegraphics[width=2.65in]{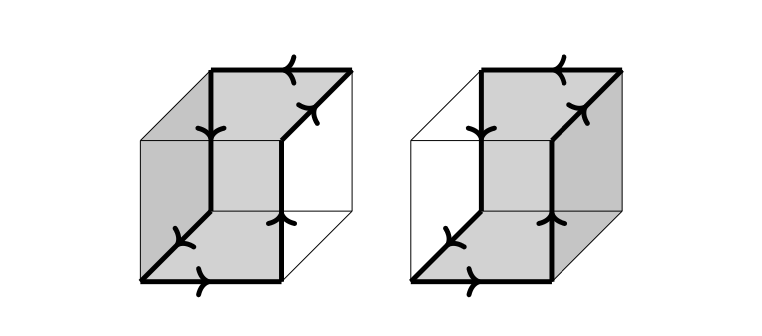}}
\caption{Two Seifert surfaces for the difference
of the two tilings of a trit move.}
\label{fig:Seifert}
\end{figure}
\end{centering}

Figure~\ref{fig:Seifert} shows two discrete Seifert surfaces for the same boundary.  The boundary is the difference of two tilings that differ in a single trit move. In the paper~\cite{FKMS}, the authors show the following key technical Lemma. 

\begin{lemma}
\label{lemma:surface}
If $\Flux(t_0) = \Flux(t_1)$ then for sufficiently large $k$ there exists a discrete Seifert surface for the pair after $k$ refinements. 
\end{lemma}

Informally, we think of refinements as giving a little more space to move;
formally, we need refinements because of Lemma~\ref{lemma:surface}.

The \emph{twist} of a tiling will prove to be a fundamental statistic for tilings, not just for the connectedness result below, see Section~\ref{sec:twist}.  The definition is below, the idea is as follows.  Given a region, fix a base tiling.  For a tiling $t$, take a Seifert surface $S$ for the pair $t$ and the base tiling.   At each vertex of $\mathcal{G}(\cR)$ we record a $0$ or $\pm 1$  depending on  how the tiles of $t$ intersect with the surface $S$.  There is an additional sign recorded, denoted $c(v)$ below, based on the color of the site.  Finally, we take the sum of all values across all vertices.

Define:
\[ \phi(t;S) = \sum_v \varphi(v;t;S) \, \, \, \, \qquad
  \varphi(v;t;S) = \textrm{c}(v) \cdot
\begin{cases}
+1, &\text{end above } S \\
0, &\text{end on } S \\
-1, &\text{end below } S
\end{cases} \, \, \, \, \qquad
\]

\begin{defin}
\label{def:twist}
Fix a base tiling $t_{\oplus}$.  The \emph{twist} of a tiling $t$ is
$$\Twist(t) := \phi(t; t-t_{\oplus}) \, \, $$
\end{defin}

The flux and twist behave nicely with respect to our local moves.  
If $t_0$ is connected to $t_1$ by a sequence of flip and trit moves:
$$ t_0 \leadsto \textrm{trit} \leadsto t_1 $$ then the flux remains the same: 
$$\Flux(t_0) = \Flux(t_1).$$
{If the sequence includes exactly one trit then}
the twist changes by exactly $1$:
$$\Twist(t_0) =  \Twist(t_1) \pm 1;$$
{the sign is defined by the orientation of the trit.}

The main theorem of~\cite[Theorem 1]{FKMS}
shows how tilings are connected by flips and trits
with respect to  flux and twist. 

\goodbreak

\begin{theo}
  \label{thm:main} 
For two tilings $t_0$ and $t_1$ of a region $\cR$:\\
\begin{enumerate}
\item{There exists a sequence of flips and trits
connecting refinements of $t_0$ and $t_1$
\emph{if and only if}
$\Flux(t_0) = \Flux(t_1)$.}
\item{There exists a sequence of flips 
connecting refinements of $t_0$ and $t_1$
\emph{if and only if}
$\Flux(t_0) = \Flux(t_1)$ and 
$\Twist(t_0) = \Twist(t_1)$.  }
\end{enumerate}
\end{theo}

The theorem shows that tilings can be connected by flips and trits if
we allow for
{(a sufficiently large number of)} refinements. 
One might think that the refinements make connectivity almost too easy,
namely that with enough refinements everything is connected. 
Theorem~\ref{thm:main} shows this is not the case. 
If two tilings have different values of the twist,
then the two tilings can not be flip connected
regardless of how many refinements are made.

 Missing from this result is how many refinements are needed to find a sequence of local moves when such a sequence exists.  Although the authors believe the number to be quite small, it remains an open problem to even bound the number of refinements needed.   

 \subsection{A Second Topological Approach}
 \label{sec:localcomplex}
 
 Our second approach
is also topological in nature. 
This approach uses cell complexes, homotopy theory and fundamental groups.

 As in the last section, we provide results given a little ``extra space''. 
Previously we used refinements,
here we  consider adding extra vertical space to cylinders. 
A \emph{cylinder} is a region $ \cR = \mathcal{D} \times [0,N]$ where
$\mathcal{D}$ is a fixed two dimensional region.
We assume $\mathcal{D}$ to be connected and simply connected,
a topological disk.
Extra vertical space is added by moving from
$\cR = \mathcal{D} \times [0,N]$ to
$\cR' = \mathcal{D} \times [0,N+M]$. 
A tiling of $\cR$ is given space
by adding all vertical dominos in the top $M$ floors.

We begin by defing the \emph{domino complex} $\mathcal{C_{\D}}$ of a
disk $\D$. $\mathcal{C_{\D}}$ is a $2$-dimensional CW-complex
constructed such that tilings of $\cR = \D \times [0,N]$ are closed
paths in the complex.  The vertices, edges, and two-cells correspond
to tilings as follows.  First, when considering a floor of a tiling,
some tiles lie completely within the floor while others intersect
vertically occupying only one square of the floor.  Vertices
correspond to the partial tilings of a floor marking only
vertical intersections.
{More formally, a \textit{plug} 
a balanced set of squares in $\D$ and there is one vertex
in our CW-complex for each plug of $\D$.}
There is an edge between two 
{disjoint plugs}
for each way
{the union}
can be completed to a full tiling of the floor. 
In this way, edges correspond to sets of
{disjoint}
horizontal tiles that lie within a floor.

\begin{figure}[ht]
    \centerline{\includegraphics[width=4.3in]{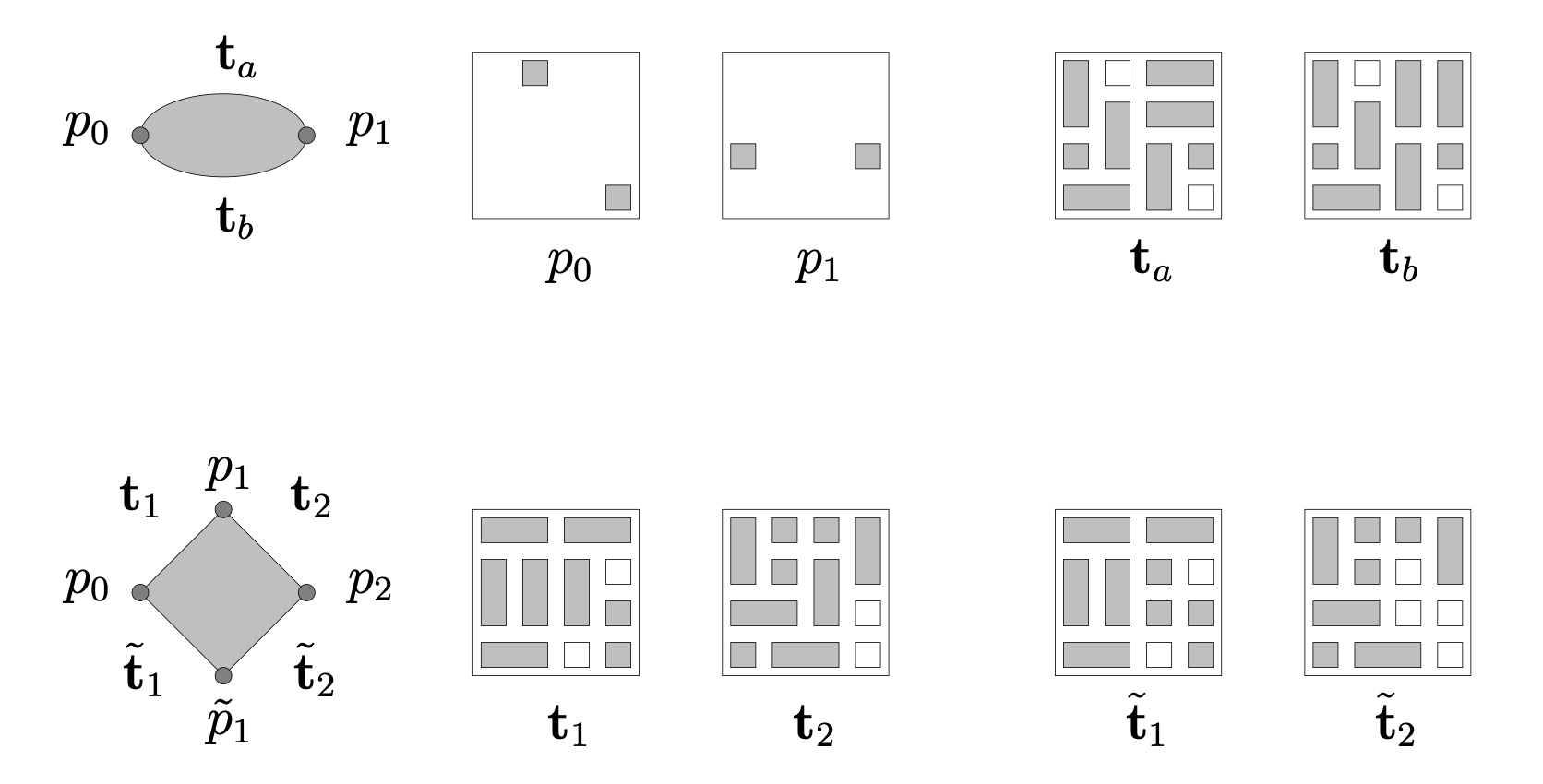}}
    \caption{The two types of $2$-cells of the domino complex.}
   \label{fig:Dcomplex}
 \end{figure}

There are two types of $2$-cells both of which reflect flip moves. 
Glue a bigon for each flip within a floor
and glue a quadrilateral for each vertical flip move,
see Figure~\ref{fig:Dcomplex}.

 The importance of the construction is in the following result.
 \begin{prop}
   Two tilings of a cylinder
are flip connected
(with possible vertical space added)
if and only if their paths in the domino complex  are homotopic.
 \end{prop}

\begin{defin}
The \emph{domino group}  $\mathcal{G_D}$ is the fundamental group of the domino complex:
$$\mathcal{G_D} = \pi_1(\mathcal{C_D}).$$
\end{defin}

The twist of a tiling defined a group homomorphism
from the domino group to the integers.
Define the homomorphism $\Tw$ by appending the parity of the height $N$
of the region:
$$\Tw: \mathcal{G_D} \rightarrow \mathbf{Z} \oplus \mathbf{Z}/2.$$
This homomorphism gives us a way to detect whether tilings are flip connected. 

\begin{prop}
The homomorphism $\Tw$ is an isomorphism if and only if
any two tilings of the same twist
are flip connected (with possible vertical space added).
\end{prop}

In higher dimension, twist assumes values in $\mathbf{Z}/2$
and $\Tw$ is a group homomorphism of the form
$\Tw: \mathcal{G_D} \rightarrow \mathbf{Z}/2 \oplus \mathbf{Z}/2$.
It has been shown that the twist is an isomorphism
for a number of different regions as collected below,
see~\cite{MR4448575, 4Dtiles, Raphael}.

\begin{theo}
\label{thm:complexes}
Twist is an isomorphism for:
\begin{enumerate}
\item Two dimensional boxes with sides at least $4$;
\item Arbitrary dimensional boxes (except a few very small examples);
\item Hamiltonian disks without bottlenecks.
\end{enumerate}
\end{theo}

\begin{figure}[ht]
    \centerline{
    \includegraphics[width=3.75in]{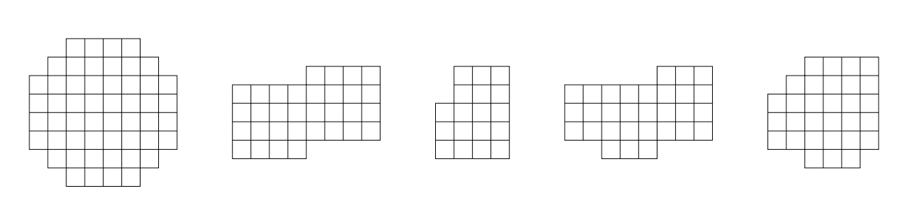}}
    \caption{Various disks for which the Twist is an isomorphism.} 
\label{fig:regular}
\end{figure}

Figure~\ref{fig:regular} shows a number of disks $\D$
such that the Twist is an isomorphism. 
Therefore, for any of the cases listed in Theorem~\ref{thm:complexes}
or shown in Figure~\ref{fig:regular},
any two tilings of the same twist of a cylinder above the region
are flip connected,
provided we are allowed to add a little extra vertical space. 

It remains an open question to understand generally for which regions
Theorem~\ref{thm:complexes} holds. 
It seems to be true for regions without small pinch points, or bottlenecks. 
However current proofs rely on an additional niceness assumption
such as having a Hamiltonian path or cycle  
(see~\cite{Raphael} for a precise statement,
including a definition of bottleneck).  
Notice that the existence of a Hamiltonian path or cycle
is sometimes a delicate question.

The amount $M$ of extra vertical space required is a function of
$\mathcal D$ only,
that is, it does not depend on the height $N$ of the region.
The proof of this fact does not give a practical way to compute $M$.
It would of course be interesting to obtain explicit estimates,
at least for such examples as the square ${\mathcal D} = [0,4]^2$.
Experimental evidence seems to indicate that a small value of $M$
would already work.

For instance, a brute force computation verifies the following fact.
For a tiling $\bt$ of ${\mathcal R} = [0,4]^2 \times [0,4]$,
let $\bt'$ be the tiling of ${\mathcal R}' = [0,4]^2 \times [0,6]$
obtained from $\bt$ by adding two floors of vertical dominoes.
Notice that $\Twist(\bt') = \Twist(\bt)$.
If $\bt_0, \bt_1$ are tilings of $\mathcal R$
and $\Twist(\bt_0) = \Twist(\bt_1)$
then there is a finite sequence of flips
taking $\bt_0'$ to $\bt_1'$ (in ${\mathcal R}'$).

\subsection{A Commutative Algebra Approach}

In this section we describe a different approach to higher dimensional
tilings involving techniques from commutative algebra, as introduced in
{Chin~\cite{Chin}, Gross, Yazmon~\cite{MR4294207}}.
The central idea to this approach is
that tilings are encoded as monomials in a polynomial ring.  Local
moves between tilings are then interpreted by constructing binomials.
As is typical in commutative algebraic approaches in combinatorics,
one forms ideals out of the objects of interest.  In this way, the
rich literature of monomial and binomial ideals can be applied to
questions in domino tilings. 
We construct two
such ideals -- the tiling ideal and the flip ideal.

Let $\cR \subset \R^N$ be a balanced cubical region and let $\GR = \mathcal{G}(R^*)$ be the graph dual of $\cR$.  In the commutative algebra setting, we will abuse notation slightly 
 and consider domino tilings of $\cR$ in terms of the corresponding
 perfect matching in $\GR$.  Associate a variable $x_e$ to each edge $e \in E(\GR)$ and consider the polynomial ring generated by
 the edge variables:
$$k[E] = k[x_e \, | \,  e \in E(G)].$$

 Let $T$ be a tiling of the region $\cR$.  The monomial of $T$ is the product of all edges in the tiling:

 $${\bf{x}}^T = \prod_{x_e \in T} x_e $$
 
 Let $T_1$ and $T_2$ be two tilings of a region $\cR$.  Define the
 difference of $T_1$ and $T_2$ as the binomial:

$$x^{T_1} - x^{T_2} = \prod_{x_e \in T_1} x_e - \prod_{x'_{e} \in T_2} x'_{e}$$ 

 As above, the difference of two tilings is a system of cycles in $\GR$.  

 \begin{defin}
   The tiling ideal of $\cR$ is the binomial ideal
   $$I_{\textrm{Tiling}} :=  \langle x^{T_1} - x^{T_2} \, | \, T_1, T_2 \textrm{ are tilings of } \cR  \rangle$$
\end{defin}

 \begin{defin} The flip ideal of $\cR$ is the ideal
   $$I_{\textrm{Flip}} := \langle x_ax_b - x_cx_d  \,  | \, (ab, cd) \textrm{ is a flip move of } \cR \rangle$$
   where $(ab,cd)$ is a flip move if the edges $a,b,c,d$ form a $4$-cycle in $\GR$ with edges $a,b$ opposite each other and edges $c,d$ opposite each other. 
   \end{defin}

   \begin{figure}[ht]
     \centerline{
\includegraphics[width=1.25in]{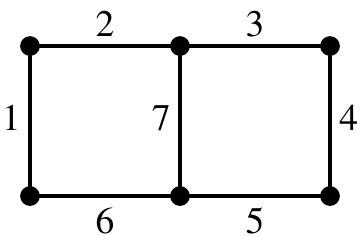}}
\caption{A small section of $\GR$ with edge labels.} 
     \label{fig:ideals}
   \end{figure}

 \begin{example}
   \label{ex:ideals}
   Figure~\ref{fig:ideals} shows two adjacent $4$-cycles which is $\GR$ for the $2 \times 3$ box.  The term associated to each $4$-cycle  is an element of the flip ideal.  For example,   $$x_4x_7 - x_3x_5 \in I_{\textrm{flip}}$$ reflecting that dominos in positions corresponding to edges $4$ and $7$ could be flipped to positions $3$ and $5$.

The tiling ideal of the region in Figure~\ref{fig:ideals} contains the term $$x_1(x_4x_7 - x_3x_5) = x_1x_4x_7 - x_1x_3x_5 \in I_{\textrm{tiling}}$$ reflecting two tilings that differ only by the flip across the cycle on the right. 
   
   \end{example}

 The connectivity of a region by flip moves can be determined by comparing the two ideals.  As in the example, one always has that a flip binomial can be extended to an element of the tiling ideal by multiplying by all other variables held fixed by the flip move.    We now state the main theorem of this approach, Theorem 3.12~\cite{MR4294207}.

 \begin{theo}
   \label{thm:binomials}
   A region $\cR$ is flip connected if and only if
$$I_{\textrm{Tiling}} \subseteq I_{\textrm{Flip}}. $$ 
 \end{theo}

  Using Theorem~\ref{thm:binomials}, the authors are able to recover the connectivity result in two dimensions:

 \begin{theo} 
  Let $\cR$ be a 2-dimensional simply connected cubiculated region. Any
  binomial arising from two distinct tilings of $\cR$ is generated by quadratics. In particular, $I_{\textrm{Tiling}} \subseteq I_{\textrm{Flip}}.$
  \end{theo}

 The work~\cite{Chin} explores primary decompositions of the flip ideal.  Informally, a primary decomposition breaks  an ideal into its irreducible components.  The ideal is then perhaps better understood by understanding its building blocks.  

 A complete primary decomposition is given for the flip ideal of
 narrow grids and a conjecture is given for the general grid case.  The idea is that the prime components are reflected by neighborhoods of internal vertices and boundary conditions.

\section{Random Tilings}

What does a random tiling look like?
As we have already mentioned above,
this is a very active area of research concerning domino tilings
in dimension two.
As an example,
Cohn, Kenyon and Propp showed
a variational principle for domino tilings in the plane~\cite{variational}.     

There is also some very interesting current work in higher dimensions.
The approach and techniques are sufficiently different from the remainder of the survey that we can not really do justice to this area of research here.
We simply mention 
a recent achievement in the probability theory of higher dimensional tilings:
a large deviation result by Chandgotia, Sheffield and Wolfram~\cite{Large}. 

  \begin{figure}[ht]
    \begin{center}
\includegraphics[height=2.4in]{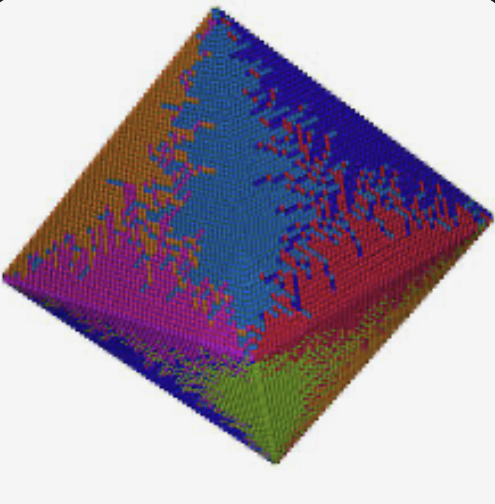}
\caption{A random tiling of the Aztec octahedron}
\label{fig:arcticocta}
\end{center}
  \end{figure}

Figure~\ref{fig:arcticocta} shows a random tiling of
the ``arctic octahedron''.    
The traditional algorithms used to generate such figures
involve random walks in the space of tilings in order to uniformly sample a tiling.
Such a random walk may use either local moves
or cycles of unbounded length.
Using local moves leads naturally to conjectures 
such as Question~\ref{qn:localbox}.

\section{Twist}
\label{sec:twist}

As we have seen already above and will be further motivated below, the twist is an important statistic for higher dimensional tilings. We note that although defined topologically, for sufficiently nice regions, the twist can be computed purely combinatorially.  This allows us to gather some data on the twist.

\begin{figure}[ht]
  \label{fig:POSET}
\begin{center}
\includegraphics[height=3in]{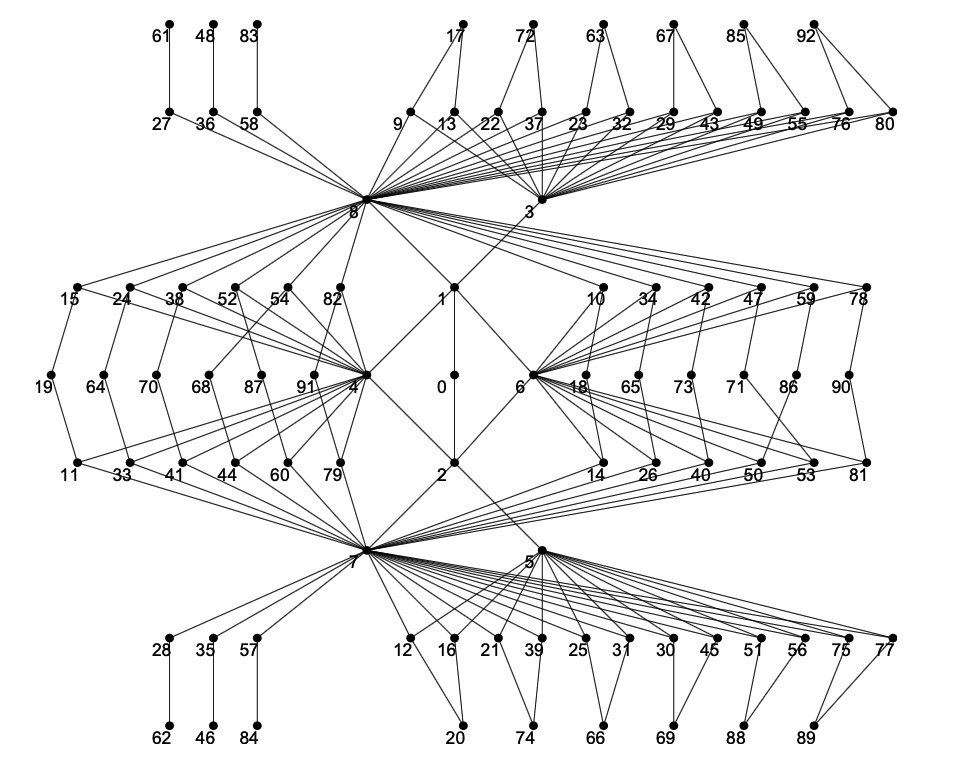}
\caption{The graph of connected components under flips
of the $4\times 4\times 4$ box.}
\end{center}
\end{figure}

Figure~\ref{fig:POSET} shows the graph of connected components under the flip move for tilings of the $4 \times 4 \times 4$ box.
The $4 \times 4 \times 4$ box has $93$ flip connected components, they
are listed with multiplicity in the table of Section~\ref{sec:basics}. The graph of
Figure~\ref{fig:POSET} has one vertex per connected component.  The
largest component is labeled $0$.  Two flip connected components $A$
and $B$ are connected by an edge in the graph if there exists tilings
$T_1 \in A$ and $T_2 \in B$ such that $T_1$ and $T_2$ are connected by
a trit move. Note there may be many such tilings, but we only draw one
edge.

The graph is laid out so that vertices in the same row have the
same value of twist.  Vertices corresponding to tilings of twist $0$
are across the middle. Tilings one level up have twist $1$.  Tilings
one level down have twist $-1$.  Recall that a trit move changes the
twist by exactly $1$, hence we only have edges between adjacent rows.
The vertical symmetry of the graph is expected due to the symmetry of
the region.  Other aspects of the graph, such as the varying number of
components per level or the varying degrees of components is not
easily anticipated.

In this example, the maximum and minimum values for the twist are $\pm 4$.  
For a general box, the possible values of the twist can be bounded as
follows, see~\cite[Proposition 6.1]{Milet}.

\begin{prop}

  For $\cR = [0,L] \times [0,M] \times [0,N]$,
  $$ \frac{1}{162} \leq \frac{\Tw_{\max}(\cR)}{LMN\min(L,M,N)} \leq \frac{1}{16}$$ 
  
  \end{prop}

\section{Helicity}

Let $M \subset \R^3$ be a simply connected domain
with smooth boundary $\partial M$.
Let $\xi$ be a divergence-free vector field in $M$
which is tangent to its boundary $\partial M$.
The \textit{helicity} of $\xi$ is defined by
\[ \Hel(\xi) = \int_M (\xi, \curl^{-1}\xi) dxdydz. \]
Thus, helicity is an invariant and it can be interpreted
as measuring the linkage and knottedness of orbits.
More precisely, in \cite{Arnold73} (see also \cite{AK}),
Arnold proved that an asymptotic version of
the Hopf invariant for a vector field is equal to the field's helicity.

If $\xi$ is not tangent to the boundary,
we can still define the \textit{relative helicity} of $\xi$.
In a nutshell, extend $\xi$ to a thin tubular neighborhood of $\partial M$
to define an \textit{isolating shell}.
Among vector fields with the same behavior at the boundary,
helicity is then well defined up to an additive constant.

\begin{figure}[ht]
\centering
\includegraphics[width=3in]{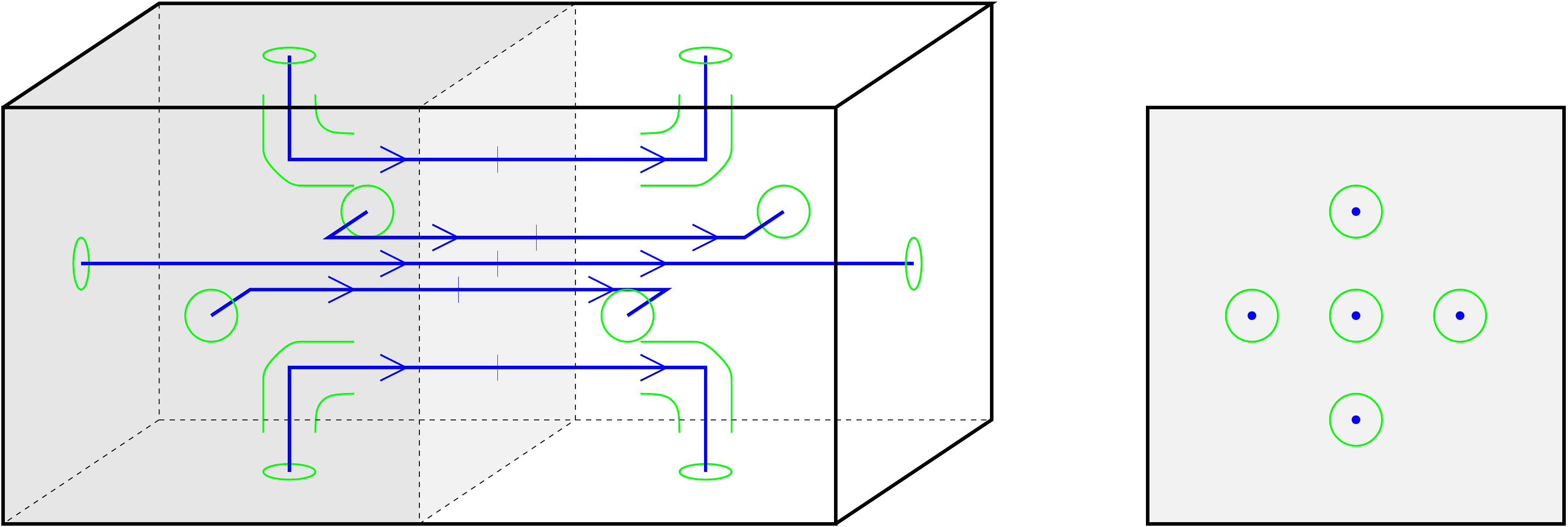}
\caption{The $5$ pipe construction (one domino).}
\label{fig:5pipe}
\end{figure}

The similarity between the concepts of twist (for tilings)
and helicity (for vector fields)
has been informally remarked several times,
including by a few audience members.
An explicit construction relating the two concepts
is the \textit{five pipe construction}, see~\cite{KhS}.
Consider a domino tiling of a 3d region $\cR$.
Inside each domino construct five narrow pipes
joining faces of the domino as in Figure~\ref{fig:5pipe}.

Notice that in the interior of $\cR$ pipes connect nicely.
At the boundary $\partial \cR$,
the pipe pattern is independent from the tiling.
Turn the water on, with flux $\varphi = 1/6$ in each pipe.
The water defines a vector field $\xi$,
which is clearly divergence-free.
We claim that the relative helicity agrees with the twist which are both are defined up to additive constant:
this is the main result (Theorem~1.1) from \cite{KhS},
which we restate in a slightly simplified form.

\begin{theo}
Let $\cR \subset \R^3$ be a cubiculated compact region.
We assume that $\cR$ is tileable and simply connected. 
Let $\bt_0$ and $\bt_1$ be tilings of $\cR$:
perform the $5$ pipe construction with flux $\varphi = 1/6$
to obtain vector fields
$\xi_{\bt_0}$ and $\xi_{\bt_1}$. We then have
\[ \Hel(\xi_{\bt_1}) - \Hel(\xi_{\bt_0}) = \Tw(\bt_1) - \Tw(\bt_0). \]
\end{theo} 

\begin{figure}[ht]
\centering
\includegraphics[width=3.5in]{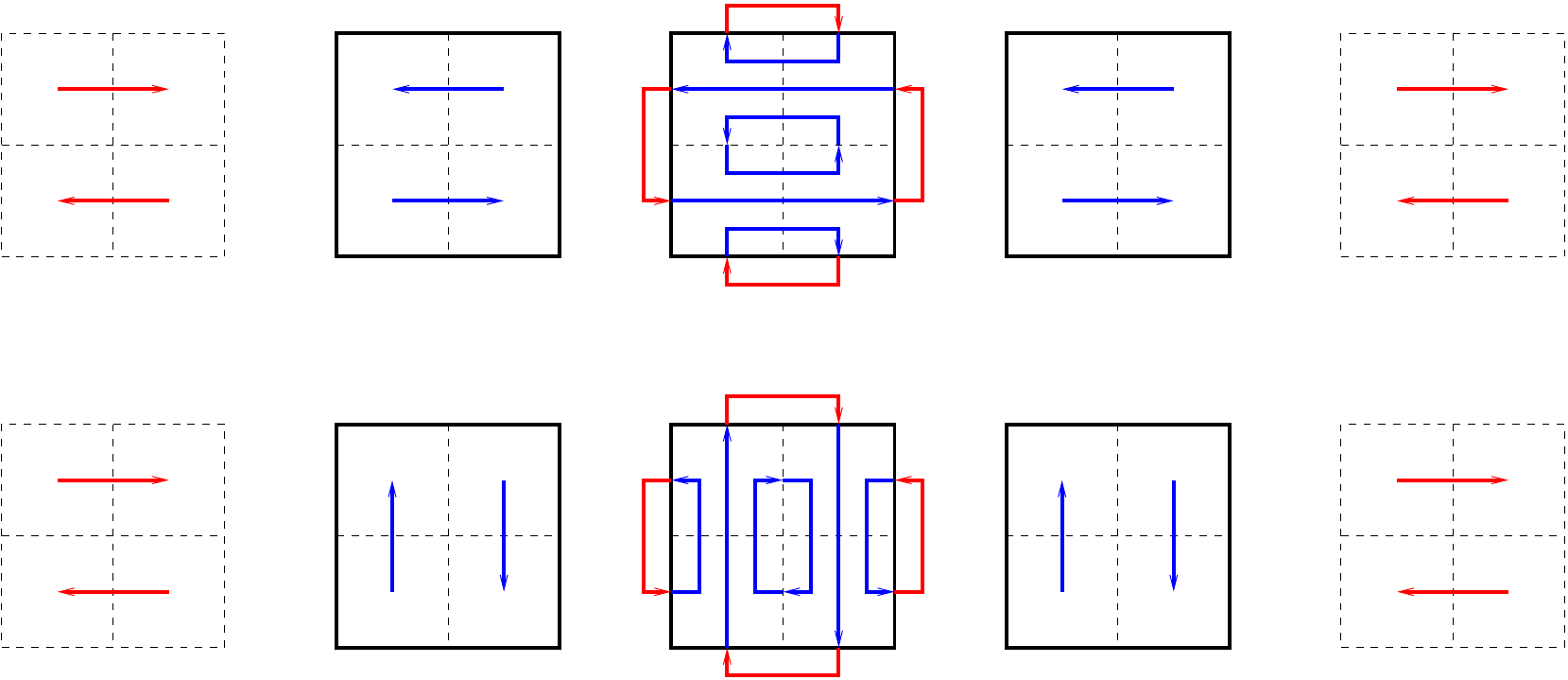}
\caption{The pipe systems corresponding to a flip.}
\label{fig:flip5c2}
\end{figure}

The proof of the theorem relies on the results of~\cite{FKMS}.
More precisely, 
we need to compute the change of relative helicity
when we perform a flip, a trit or a refinement.
After many such figures we verify that, as expected,
relative helicity remains unchanged under flip or refinement
and changes by $\pm 1$ under a trit (with correct sign).
As an example, we show in Figure~\ref{fig:flip5c2}
the pipe systems and vector fields corresponding to a flip:
the relative helicity does not change.

\section{Enumeration}
\label{sec:count}

A natural question is: how many domino tilings are there of a given
region in $3$-dimensions?  The analogous question in $2$-dimensions
has been studied extensively both for domino and other types of
tilings.  Enumerating such tilings proves to be quite difficult and
leads to fascinating connections and results. The Chapter by Krattenthaler of the
current volume investigates the enumeration questions of tiling in
$2$-dimensions in great depth.

As one might imagine, there is less known about the case in higher
dimensions.  Recall the table of Section~\ref{sec:basics} giving
the sizes of the flip connected components for small boxes.
\begin{itemize}
\item[]The
  number of tilings of the $3 \times 3 \times 2$ box is $229$.
\item[]The number of tilings of the $4 \times 4 \times 4$ box is
$5,051,532,105$.
\item[]The number of tilings of the $4 \times 4 \times 8$ box is
\[ 175,220,727,982,196,365,632. \]
  \end{itemize}

There are a tremendous number of tilings in higher dimensions, growing rapidly with the size of the box.  In general, we do not have a formula for the number of tilings of a given region, even for a box.
\begin{question}
  \label{qn:howmany}
 How many tilings are there of the d-dimensional box?
\end{question}

Unfortunately, due to the next result, it is unlikely that a nice
formula can be found for the number of domino tilings in higher
dimensions.  First proved by Valiant~\cite{Valiant} and later strengthened by Pak and Yang~\cite{PakYang} is the following complexity result.

\begin{theo}
The number of domino tilings of a region in $3$-dimensions is \#P-complete. 
\end{theo}

Although a perfect count may be out of reach, understanding  the
approximate number of domino tilings is an interesting question and important in applied physical
contexts.  A number of works have proved upper and lower bounds on the number of tilings of the $d$-dimensional box.  We will not discuss this line of inquiry here but refer the reader to~\cite{MR1635888, MR237350, MR560565}.

In two dimensions, many of the classical enumerative techniques are
linear-algebraic in nature.  For example, one constructs a signed
weighted adjacency matrix known as the Kasteleyn matrix and then takes
a determinant in order to enumerate tilings of rectangular regions.  Such a construction can be mimicked for $3$-dimensional boxes.  The
result is not an exact count but it is an interesting statistic on
tilings.  For a box $\cR$, let $K(\cR)$ be the generalized Kasteleyn matrix, which is essentially a signed adjacency matrix with entry $K_{ij} = 0$ if no domino occupies the cubes $i$ and $j$ and is $\pm 1$ depending on the color of the sites $i$ and $j$.    

  \begin{theo}
    \label{thm:alternating}
  The Pfaffian of the matrix $K(\cR)$ gives an alternating sum of
  tilings by the parity of the twist.  $$ \pf({K(\cR)}) = \sum_{T} (-1)^{\Tw(T)} $$  
\end{theo}

  Theorem~\ref{thm:alternating} reinforces a natural question refining
  Question~\ref{qn:howmany}: How many tilings are there of the
  $d$-dimensional box for a fixed value of the twist? Given the difficulty of exact counting and the large number of tilings,   we rescope the problem as understanding the \emph{distribution} of twist for a random tiling of the box.  

  Again looking at the data of the table in Section~\ref{sec:basics}, one might be
  tempted to conjecture that there is a giant component in the graph
  of flip connected components and that almost all tilings of the box
  have twist equal to $0$.  However, we know that this is not the
  case.

\begin{figure}

   \begin{center}
\includegraphics[height=1.5in]{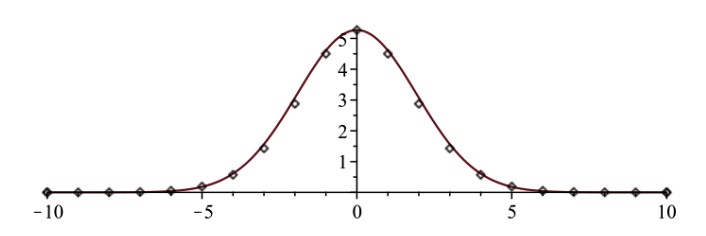}
\caption{$\#$Tilings per value of Twist of the $4 \times 4 \times 60$ box.}
\label{fig:normal}
   \end{center}
 \end{figure}

  Figure~\ref{fig:normal} plots the number of tilings by twist
  of the $4 \times 4 \times 60$ box.  The vertical axis is scaled by a
  factor of $10^{156}$ and the curve drawn in red is the true
  Gaussian.  
 For cylinders, i.e. regions of the form $\cR = D \times [0,N]$
where $D$ is a fixed two dimensional disk,
it is known that the twist is normally distributed as $N \rightarrow \infty$,
see~\cite[Theorem $3$]{normal}.  
The distribution of the twist is unknown for arbitrary boxes
when more than one dimension goes to infinity.

\begin{question}\label{qn:normalLMN}
Is the twist  normally distributed for boxes
$\R = [0,N] \times [0,M] \times [0,L]$ as $\min{N,M,L} \rightarrow \infty$ ? 
\end{question}

\section{Slabs}
\label{sec:slabs}

\begin{figure}[ht]
\begin{center}
\includegraphics[width=4.75in]{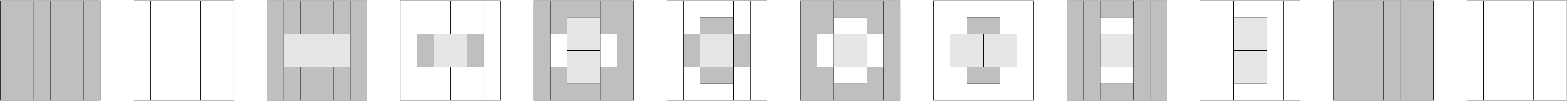}
\caption{A slab tiling of a box.}
\label{fig:redgreen4a}
\end{center}
\end{figure}

Very little is known about tilings by other higher dimensional tiles
but there are a few results about tilings of boxes by \textit{slabs},
$2 \times 2 \times 1$ blocks.
We present a small sample of such results, see~\cite{DCG}.
Figure \ref{fig:redgreen4a} shows an example of a slab tiling.

\begin{figure}[ht]
\begin{center}
\includegraphics[width=4.75in]{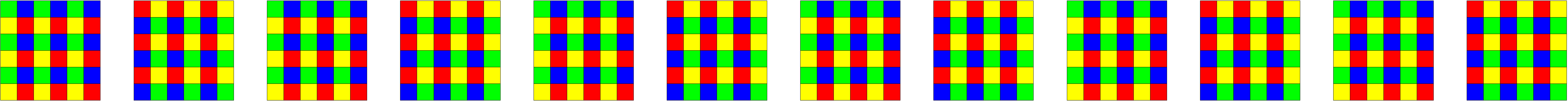}
\caption{Painting unit cubes using four colors}
\label{fig:redgreen4b}
\end{center}
\end{figure}

A \textit{flip} is a local move involving two slabs.
Is the tiling in Figure~\ref{fig:redgreen4a} equivalent
via a finite sequence of flips to a horizontal tiling?
No! There exists an invariant, the \textit{triple twist},
assuming values in ${\Z}^3$.
We proceed to compute one of its coordinates
for the example in Figure~\ref{fig:redgreen4a}.
Begin by painting the unit cubes in four colors,
red, yellow, green and blue,
as in Figure~\ref{fig:redgreen4b}.

\begin{figure}[ht]
\begin{center}
\includegraphics[width=4.75in]{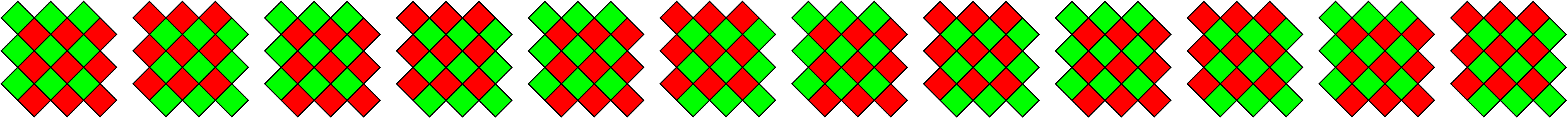}
\caption{Inflate red and green cubes}
\label{fig:redgreen4c}
\end{center}
\end{figure}

Choose two colors, say red and green.
Inflate the red and green cubes; deflate the yellow and blue. 
The result is shown in Figure~\ref{fig:redgreen4c}.

\begin{figure}[h!]
\begin{center}
\includegraphics[width=4.75in]{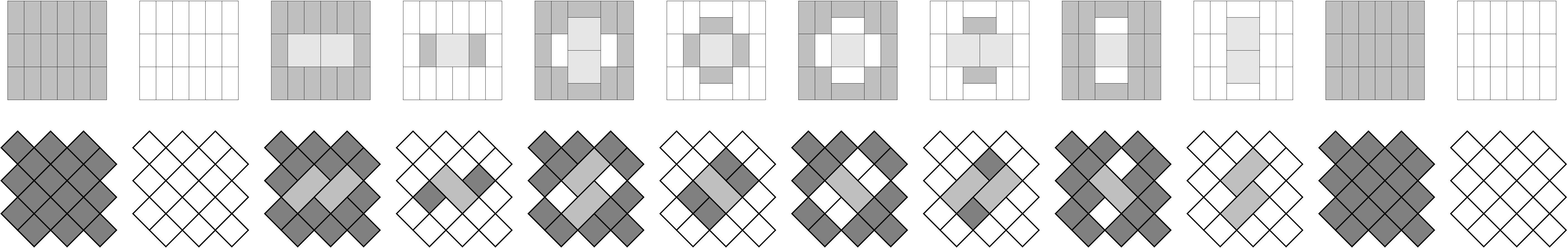}
\caption{The slab tiling becomes a domino tiling}
\label{fig:redgreen4d}
\end{center}
\end{figure}

The slab tiling in Figure~\ref{fig:redgreen4a}
becomes a domino tiling,
as in Figure~\ref{fig:redgreen4d}.

Compute the  twist: the result is invariant by flips.
Notice that the horizontal tiling yields twist zero
while the example in the figures yields twist $\pm 2$,
proving the claim above.  Perform the same construction 
for other pairs of colors to define a total of six 
related invariants.
They turn out not to be linearly independent,
but we still have three dimensions, thus defining the triple twist.

\goodbreak

\section{Four dimensions and above}

The two dimensional theory of domino tilings is well studied and we have seen above how many ideas can be extended to the
three dimensional case, although with far more open questions than
answers.
There is a distinct shift when we move to four dimensions and above.   Borrowing from Section~\ref{sec:local}, for example, we might expect a new topological invariant and corresponding local move as we go up in dimension.  This is, however, not the case.

In dimensions $n \geq 4$, a domino is a $2 \times 1 \times \cdots \times 1$ block. In the paper~\cite{4Dtiles}, the authors establish the following theory.   The \emph{twist} of a tiling is defined; it is no longer an integer but an element of $\Z/(2)$.  On the one hand this is simpler than the three dimensional case, on the other hand it collapses much of the interesting behavior seen above.  The \emph{flip} and \emph{trit} moves remain the same, involving either two or three dominos and take place in a two or three dimensional subspace of the tiling region.  The flip move does not alter the value of the twist, a trit move changes the value of the twist  to the only other possible value.

In dimensions $n \geq 4$, there is both a connectedness-by-local-moves theorem and an enumerative theorem.  The connectedness result is stated in Section~\ref{sec:localcomplex}.  In particular, the result uses the theory of domino complexes and the domino group for cylinders $\cR_{\mathcal{D}} = \mathcal{D} \times [0,N]$ where $\mathcal{D}$ is a region in $\R^{n-1}$.   Reiterated here,
\begin{theo}
  For boxes $\mathcal{D} = [0,N_1] \times \cdots \times [0,N_{n-1}]$, any two tilings of $\cR_{\mathcal{D}}$ with the same twist can be connected by a sequence of flip moves (with possible vertical space added).
\end{theo}

In this case, the amount of extra vertical space
needed can be bounded in terms of the disk $\mathcal{D}$. 
These results then allow
us to understand the distribution of the twist in higher dimensions. 
We give a simplified version here,
see~\cite[Corollary 1.2]{4Dtiles} for the full statement. 

\begin{theo}
For boxes $\mathcal{D}$,
the set of flip connected components of tilings
of the cylinder $\cR_{\mathcal{D}}$ has two twin giant components,
one for each value of the twist.
{The ratio between the size of the two components tends to $1$
when the height $N$ goes to infinity.}
\end{theo}

\bibliographystyle{plain}
\bibliography{refs}

\begin{thebibliography}{10}

\bibitem{DCG}
Vieira Arthur M.~M. Alencar George L.~D., Saldanha Nicolau~C.
\newblock Slab tilings, flips and the triple twist.
\newblock {\em Discrete \& Computational Geometry}, 2025.

\bibitem{Arnold73}
V.~I. Arnold.
\newblock The asymptotic {H}opf invariant and its applications.
\newblock volume~5, pages 327--345. 1986.
\newblock Selected translations.

\bibitem{AK}
Vladimir~I. Arnold and Boris~A. Khesin.
\newblock {\em Topological methods in hydrodynamics}, volume 125 of {\em
  Applied Mathematical Sciences}.
\newblock Springer, Cham, second edition, [2021] \copyright 2021.

\bibitem{Large}
Sheffield Chandgotia and Wolfram.
\newblock Large deviations for the 3d dimer model.
\newblock {\em Arxiv Preprint arXiv:math/230408468}.

\bibitem{Chin}
Tracy Chin.
\newblock A computational commutative algebra approach to tilings.
\newblock Bachelor's thesis, Brown University, 2019.

\bibitem{MR1635888}
Mihai Ciucu.
\newblock An improved upper bound for the {$3$}-dimensional dimer problem.
\newblock {\em Duke Math. J.}, 94(1):1--11, 1998.

\bibitem{Aztec2}
Henry Cohn, Noam Elkies, and James Propp.
\newblock Local statistics for random domino tilings of the {A}ztec diamond.
\newblock {\em Duke Math. J.}, 85(1):117--166, 1996.

\bibitem{variational}
Henry Cohn, Richard Kenyon, and James Propp.
\newblock A variational principle for domino tilings.
\newblock {\em J. Amer. Math. Soc.}, 14(2):297--346, 2001.

\bibitem{Raphael}
Raphael de~Marreiros.
\newblock Domino tilings of $3$-dimensional cylinders.
\newblock Phd thesis, Pontificia Universidade Catolica do Rio de Janeiro, 2025.

\bibitem{FKMS}
Juliana Freire, Caroline~J. Klivans, Pedro~H. Milet, and Nicolau~C. Saldanha.
\newblock On the connectivity of spaces of three-dimensional domino tilings.
\newblock {\em Trans. Amer. Math. Soc.}, 375(3):1579--1605, 2022.

\bibitem{MR4294207}
Elizabeth Gross and Nicole Yamzon.
\newblock Binomial ideals of domino tilings.
\newblock {\em Discrete Math.}, 344(11):Paper No. 112530, 14, 2021.

\bibitem{MR237350}
J.~M. Hammersley.
\newblock An improved lower bound for the multidimensional dimer problem.
\newblock {\em Proc. Cambridge Philos. Soc.}, 64:455--463, 1968.

\bibitem{Aztec1}
Propp Jockusch and Shor.
\newblock Random domino tilings and the arctic circle theorem.
\newblock {\em Arxiv Preprint arXiv:math/9801068}.

\bibitem{KASTELEYN19611209}
P.W. Kasteleyn.
\newblock The statistics of dimers on a lattice: I. the number of dimer
  arrangements on a quadratic lattice.
\newblock {\em Physica}, 27(12):1209--1225, 1961.

\bibitem{kenyon_2004}
R.~Kenyon.
\newblock An introduction to the dimer model.
\newblock 17:p. 267–304, Mar 2004.

\bibitem{Amoebae}
Richard Kenyon, Andrei Okounkov, and Scott Sheffield.
\newblock Dimers and amoebae.
\newblock {\em Ann. of Math. (2)}, 163(3):1019--1056, 2006.

\bibitem{KhS}
B.~Khesin and N.~Saldanha.
\newblock Relative helicity and tiling twist.
\newblock {\em Transactions of the American Mathematical Society}, 2025.

\bibitem{4Dtiles}
Caroline~J. Klivans and Nicolau~C. Saldanha.
\newblock Domino tilings and flips in dimensions 4 and higher.
\newblock {\em Algebr. Comb.}, 5(1):163--185, 2022.

\bibitem{Milet}
Pedro Milet.
\newblock Domino tilings of three-dimensional regions.
\newblock Phd thesis, Pontificia Universidade Catolica do Rio de Janeiro, 2015.

\bibitem{MR560565}
Henryk Minc.
\newblock An asymptotic solution of the multidimensional dimer problem.
\newblock {\em Linear and Multilinear Algebra}, 8(3):235--239, 1979/80.

\bibitem{PakYang}
Igor Pak and Jed Yang.
\newblock The complexity of generalized domino tilings.
\newblock {\em Electron. J. Combin.}, 20(4):Paper 12, 23, 2013.

\bibitem{normal}
Nicolau~C. Saldanha.
\newblock Domino tilings of cylinders: connected components under flips and
  normal distribution of the twist.
\newblock {\em Electron. J. Combin.}, 28(1):Paper No. 1.28, 23, 2021.

\bibitem{MR4448575}
Nicolau~C. Saldanha.
\newblock Domino tilings of cylinders: the domino group and connected
  components under flips.
\newblock {\em Indiana Univ. Math. J.}, 71(3):965--1002, 2022.

\bibitem{MR1370508}
Nicolau~C. Saldanha and Carlos Tomei.
\newblock An overview of domino and lozenge tilings.
\newblock volume~2, pages 239--252. 1995.
\newblock Combinatorics Week (Portuguese) (S\~ao Paulo, 1994).

\bibitem{MR136398}
H.~N.~V. Temperley and Michael~E. Fisher.
\newblock Dimer problem in statistical mechanics---an exact result.
\newblock {\em Philos. Mag. (8)}, 6:1061--1063, 1961.

\bibitem{MR1072815}
William~P. Thurston.
\newblock Conway's tiling groups.
\newblock {\em Amer. Math. Monthly}, 97(8):757--773, 1990.

\bibitem{Valiant}
L.~G. Valiant.
\newblock Completeness classes in algebra.
\newblock In {\em Conference {R}ecord of the {E}leventh {A}nnual {ACM}
  {S}ymposium on {T}heory of {C}omputing ({A}tlanta, {G}a., 1979)}, pages
  249--261. ACM, New York, 1979.

\end{thebibliography}

\medskip

\noindent
\footnotesize
Caroline J. Klivans \\
Brown University  \\ Box F \\ Providence, RI 02906 USA \\
\url{klivans@brown.edu}

\smallskip

\noindent
\footnotesize
Nicolau C. Saldanha\\
Departamento de Matem\'atica, PUC-Rio, \\
R. Marqu\^es de S. Vicente 255,
Rio de Janeiro, RJ 22451-900, Brazil  \\
\url{saldanha@puc-rio.br}

\end{document}